\documentclass[a4paper,11pt]{amsart}
\usepackage{amssymb,latexsym}
\usepackage{graphicx}
\newtheorem{theorem}{Theorem}
\newtheorem{lemma}{Lemma}
\newtheorem{definition}{Definition}
\newtheorem{notation}{Notation}
\newtheorem{proposition}{Proposition}

\newtheorem{algorithm}{Algorithm}

\def\supp{\mathop{\rm supp}\nolimits}
\def\hess{\mathop{\rm Hess}\nolimits}
\def\grad{\mathop{\rm grad}\nolimits}

\def\nor{\mathop{\rm nor}\nolimits}

\begin{document}

\title{Riemannian Median and Its Estimation}

\author{Le Yang}
\address{Laboratoire de Math\'{e}matiques et Applications, CNRS :
UMR 6086, Universit\'{e} de
Poitiers\\
T\'{e}l\'{e}port 2 - BP 30179, Boulevard Marie et Pierre Curie,
86962 Futuroscope Chasseneuil Cedex, France}
\email{Le.Yang@math.univ-poitiers.fr}


\begin{abstract}

In this paper, we define the geometric median of a probability
measure on a Riemannian manifold, give its characterization and a
natural condition to ensure its uniqueness. In order to calculate
the median in practical cases, we also propose a subgradient
algorithm and prove its convergence as well as estimating the error
of approximation and the rate of convergence. The convergence
property of this subgradient algorithm, which is a generalization of
the classical Weiszfeld algorithm in Euclidean spaces to the context
of Riemannian manifolds, also answers a recent question in P. T.
Fletcher et al. \cite{Fletcher3}

\end{abstract}


\thanks{The author is very grateful to his supervisors: Marc Arnaudon and
Fr\'{e}d\'{e}ric Barbaresco for many helpful discussions and their
penetrating comments. The author will also thank Thales Air Systems
for their financial supports.}

\maketitle

\section{introduction}

The classical Fermat point of a triangle in the plane is the point
minimizing the sum of distances to the three vertices. This is a
prototype of the more general Fermat-Weber problem concerning the
same question but for more than three points and in higher
dimensions. The solution of this problem is called the geometric
median of these points and provides a notion of centrality for them.
For this reason, the geometric median is a natural estimator in
statistics which possesses another important property called
robustness, that is, not sensitive to outliers. As a consequence,
the geometric median is a widely used robust estimator in both
theoretical and practical theory of robust statistics.

Naturally, one can also ask the question to find a point that
minimizes the sum of distances to a set of given points in a much
more general space as long as it carries a distance. This has been
done in Sahib \cite{Sahib} who proved the existence of the geometric
median of a probability measure on a complete, separable and
finitely compact metric space. Recently, there is a growing interest
of the method that characterize the statistical data lying on a
Riemannian manifold and its applications , see for example
Barbaresco \cite{Barbaresco1}, \cite{Barbaresco2},
\cite{Barbaresco4}, Fletcher et al. \cite{Fletcher1},
\cite{Fletcher2} and Pennec \cite{Pennec} in which the centrality of
empirical data is modeled by the Riemannian barycenter which was
first introduced by Karcher \cite{Karcher} and then has been studied
by many other authors, for example see Arnaudon \cite{Arnaudon1},
\cite{Arnaudon2}, \cite{Arnaudon3}, Emery \cite{Emery} and Kendall
\cite{Kendall}. As is known to all that the barycenter is not a
robust estimator and sensitive to outliers, in order to overcome
this drawback, Fletcher et al. \cite{Fletcher3} defined the weighted
geometric median of a set of discrete sample points lying on a
Riemannian manifold and proved its existence and uniqueness.

In many cases, especially in practice, one often needs to calculate
or at least estimate the value of the geometric median. In the case
of Euclidean spaces, Weiszfeld algorithm proposed firstly by
Weiszfeld \cite{Weiszfeld}, is a well known algorithm to do this
calculation and has been studied, improved on by many other authors,
for example see Khun \cite{Khun} and Ostresh \cite{Ostresh}. In the
contexts of Riemannnian manifolds, Fletcher et al. \cite{Fletcher3}
proposed a Riemannian generalization of Weiszfeld algorithm to
estimate their geometric median, they proved a convergence result
under the condition that the manifold is positively curved, it
should be noted that the range of their stepsize $\alpha$ should not
be $[0,2]$ but $[1,2]$, and conjectured the convergence in
negatively curved case.

The aim of this paper is to define the geometric median of a
probability measure on a complete Riemannian manifold and
investigate its uniqueness as well as its approximation. As in
Karcher \cite{Karcher} and Le \cite{Le}, we suppose that the support
of the probability measure is contained in a convex ball and we give
a characterization of the geometric median which is proved in the
case of Euclidean space by Khun \cite{Khun} for a discrete set of
sample points. Then we prove the uniqueness of the median under a
natural condition imposed on the probability measure and show that
this condition yields a strong convexity property which is useful in
error estimates. By regarding the Weiszfeld algorithm as a
subgradient procedure, we introduce a subgradient algorithm to
estimate the median and prove that this algorithm always converges
without condition of the sign of curvatures by generalizing the
fundamental inequality in Ferreira and Oliveivra \cite{Ferreira} in
which it was proved in positively curved manifolds. Finally, the
results of approximating errors and rate of convergence are also
obtained.


Throughout this paper, $M$ is a complete 
Riemannian manifold with Riemannian metric
$\langle\,\cdot\,,\cdot\,\rangle$ and Riemannian distance $d$. The
gradient operator and hessian operator on $M$ are denoted by $\grad$
and $\hess$, respectively. For every point $p$ in $M$, let $d_p$ be
the distance function to $p$ defined by $d_p(x)=d(x,p)$. We fix a
convex ball $B(a,\rho)$ in $M$ centered at $a$ with a finite radius
$\rho$, here the convexity of $B(a,\rho)$ means that for every two
points $x$ and $y$ in it, there is a unique shortest geodesic from
$x$ to $y$ in $M$ that lies in $B(a,\rho)$. The lower and upper
bounds of sectional curvatures $K$ in $\bar{B}(a,\rho)$ are denoted
by $\delta$ and $\Delta$ respectively. Since $\rho$ is finite,
$\delta$ and $\Delta$ are also finite. \emph{If $\Delta>0$, we
assume further that $\rho<\pi/(4\sqrt{\Delta})$.} It is easy to
check that the following classical comparison theorems in Riemannian
geometry: Alexandrov's theorem, Toponogov's theorem and Hessian
comparison theorem can be all applied in $\bar{B}(a,\rho)$, thus it
is necessary to introduce some notations for model spaces that
provide us many geometric informations.

\begin{notation}
Let $\kappa$ be a real number, the model space $M^2_{\kappa}$ is
defined as follows:\\

1) if $\kappa>0$ then $M^2_{\kappa}$ is obtained from the sphere
$\mathbb{S}^2$ by multiplying the distance function by
$1/\sqrt{\kappa}$;\\

2) if $\kappa=0$ then $M^2_0$ is Euclidean space $\mathbb{E}^n$;\\

3) if $\kappa<0$ then $M^2_{\kappa}$ is obtained from the hyperbolic
space $\mathbb{H}^2$ by multiplying the distance function by
$1/\sqrt{-\kappa}$.

The distance between two points $A$ and $B$ in $M^2_{\kappa}$ will
be denoted by $\bar{d}(A,B)$.
\end{notation}

\begin{notation}
Let $\kappa$ be a real number, then we write for $t\in \mathbf{R}$,
\[
S_{\kappa}(t)=
\begin{cases}
\sin(\sqrt{\kappa}\,t)/\sqrt{\kappa}   & \text{if $\kappa>0$;}\\
t   & \text{if $\kappa=0$;}\\
\sinh(\sqrt{-\kappa}\,t)/\sqrt{-\kappa}   & \text{if $\kappa<0$;}
\end{cases}
\]
\end{notation}

We begin with two useful estimations for our purposes, which are
almost direct corollaries of the Hessian comparison theorem. Observe
that the second estimation is the core of Le \cite{Le} where it is
formulated in terms of Jacobi fields and is proved using an argument
of index lemma.

\begin{lemma}\label{hessian estimation}
Let $p\in\bar{B}(a,\rho)$ and
$\gamma:\,[\,0,b\,]\rightarrow\bar{B}(a,\rho)$ be a geodesic, then\\

i) \[\qquad \hess d_p(\dot{\gamma}(t),\dot{\gamma}(t))\geq
D(\rho,\Delta)|\,\dot{\gamma}^{\nor}(t)\,|^2\]\\
for every $t\in[\,0,b\,]$ such that $\gamma(t)\neq p$, where
$D(\rho,\Delta)=S_{\Delta}'(2\rho)/S_{\Delta}(2\rho)>0$ and
$\dot{\gamma}^{\nor}(t)$ is the normal component of
$\dot{\gamma}(t)$ with respect to the geodesic from $p$ to
$\gamma(t)$.\\

ii) \[\hess \frac{1}{2}d^{2}_p(\dot{\gamma}(t),\dot{\gamma}(t))\leq
C(\rho,\delta)|\dot{\gamma}|^2\]\\ for every $t\in[\,0,b\,]$, where
the constant $C(\rho,\delta)\geq 1$ is defined by
\[
C(\rho,\delta)=
\begin{cases}
1  & \text{if $\delta\geq0$;}\\

2\rho\sqrt{-\delta}\coth(2\rho\sqrt{-\delta}) & \text{if
$\delta<0$;}
\end{cases}
\]

\end{lemma}

\begin{proof}
Since in $\bar{B}(a,\rho)$ we have $\delta\leq K\leq\Delta$, hence
by the classical Hessian comparison theorem we get, for
$\gamma(t)\neq p$,
\[\frac{S_{\Delta}'(d(\gamma(t),p))}{S_{\Delta}(d(\gamma(t),p))}|\,\dot{\gamma}(t)^{\nor}\,|^2
\leq\hess d_p(\dot{\gamma}(t),\dot{\gamma}(t))\leq
\frac{S_{\delta}'(d(\gamma(t),p))}{S_{\delta}(d(\gamma(t),p))}|\,\dot{\gamma}(t)^{\nor}\,|^2\]
Since $S_{\Delta}'(\theta)/S_{\Delta}(\theta)$ is nonincreasing for
$\theta> 0$ if $\Delta\leq0$ and for $\theta\in (0,\pi)$ if
$\Delta>0$, so the left inequality together with $d(\gamma(t),p)\leq
2\rho$ proves the first assertion. To show the second one, let
$\dot{\gamma}(t)^{\tan}$ be the tangential component of
$\dot{\gamma}(t)$ with respect to the geodesic from $p$ to
$\gamma(t)$ then
\begin{align*}
\hess \frac{1}{2}d^{2}_p(\dot{\gamma}(t),\dot{\gamma}(t))
&=d(\gamma(t),p)\hess
d_p(\dot{\gamma}(t),\dot{\gamma}(t))+|\,\dot{\gamma}(t)^{\tan}\,|^2\\
&\leq d(\gamma(t),p)
\frac{S_{\delta}'(d(\gamma(t),p))}{S_{\delta}(d(\gamma(t),p))}|\,\dot{\gamma}(t)^{\nor}\,|^2+|\,\dot{\gamma}(t)^{\tan}\,|^2\\
&\leq\max\bigg\{d(\gamma(t),p)\frac{S_{\delta}'(d(\gamma(t),p))}{S_{\delta}(d(\gamma(t),p))},\,1\bigg\}|\dot{\gamma}|^2
\end{align*}
Observe that $\theta S'_{\delta}(\theta)/S_{\delta}(\theta)\leq 1$
for $\theta>0$ if $\delta=0$ and for $\theta\in[0,\pi)$ if
$\delta>0$, thus for the case when $\delta\geq0$ we have
\[d(\gamma(t),p)\frac{S_{\delta}'(d(\gamma(t),p))}{S_{\delta}(d(\gamma(t),p))}\leq 1\]
if $\delta<0$, then $\theta
S'_{\delta}(\theta)/S_{\delta}(\theta)\geq1$ and is nondecreasing
for $\theta\geq0$, thus we get
\[d(\gamma(t),p)\frac{S_{\delta}'(d(\gamma(t),p))}{S_{\delta}(d(\gamma(t),p))}\leq
2\rho\frac{S'_{\delta}(2\rho)}{S_{\delta}(2\rho)}=2\rho\sqrt{-\delta}\coth(2\rho\sqrt{-\delta})\]
hence the estimation holds for every $\delta\in \mathbf{R}$.
Finally, the case when $\gamma(t)=p$ is trivial and the proof is
complete.
\end{proof}

\section{definition of riemannian median}

As in Karcher \cite{Karcher}, we now consider a probability  measure
$\mu$ on $M$ whose support is contained in the open ball
$B(a,\rho)$, and define a function
\[f:\qquad \bar{B}(a,\rho) \longrightarrow\mathbf{R}_+\,,\qquad x\longmapsto\int_M d(x,p)\mu(dp)\]
This function is well defined since for every $x\in
\bar{B}(a,\rho)$,
\[f(x)=\int_{B(a,\rho)}d(x,p)\mu(dp)\leq2\rho\]
and here are some simple properties of $f$:

\begin{lemma}\label{convexity of f}
$f$ has the following properties:\\
i) $f$ is 1-Lipschitz thus continuous;\\
ii) $f$ is convex.
\end{lemma}

\begin{proof}

$i)$ For every $x,y\in \bar{B}(a,\rho)$ we have
\[|f(x)-f(y)|\leq\int_{M}|d(x,p)-d(y,p)|\mu(dp)\leq\int_{M}d(x,y)\mu(dp)=d(x,y)\]
$ii)$ Let $p\in\bar{B}(a,\rho)$ and
$\gamma:\,[\,0,1]\rightarrow\bar{B}(a,\rho)$ be a geodesic then it
suffices to show that the function $t\rightarrow d(\gamma(t),p)$ is
convex. If $p\in\gamma[\,0,1\,]$, the convexity is trivial and if
not, it follows form the first estimation in Lemma \ref{hessian
estimation}.

\end{proof}

Since $f$ is continuous on the compact set $\bar{B}(a,\rho)$, it
attains its minimum here. Now it is time to give these minimum
points of $f$ a proper name. As the definition of barycenter in
Emery and Mokobodzki \cite{Emery}, we give the following definition.

\begin{definition}
The set of all the minimum points of $f$ is called the Riemannian
median (or median) of $\mu$ and is denoted by $\mathfrak{M}_{\mu}$.
The minimal value of $f$ will be denoted by $f_*$.
\end{definition}

It is easily seen that $\mathfrak{M}_{\mu}$ is compact. Since $f$ is
convex then along every geodesic, its right and left derivatives
exist and now we calculate them for later use.

\begin{proposition}\label{right derivative of f}
Let $\gamma:\,[\,0,b\,]\rightarrow \bar{B}(a,\rho)$ be a geodesic,
then
\[\frac{d}{dt}f(\gamma(t))\big{|}_{t=t_0+}=\langle\,\dot{\gamma}(t_0),\,H(\gamma(t_0))\,\rangle+\mu\{\gamma(t_0)\}|\dot{\gamma}|,\quad t_0\in[\,0,b\,)\]
\[\frac{d}{dt}f(\gamma(t))\big{|}_{t=t_0-}=\langle\,\dot{\gamma}(t_0),\,H(\gamma(t_0))\,\rangle-\mu\{\gamma(t_0)\}|\dot{\gamma}|,\quad t_0\in(0,b\,]\]
where for $x\in \bar{B}(a,\rho)$,
\[H(x)=\int_{M\setminus\{x\}}\grad d_p(x)\mu(dp)=\int_{M\setminus\{x\}}\frac{-\exp_x^{-1}p}{d(x,p)}\mu(dp)\in T_xM\]
is well defined and satisfying $|H(x)|\leq 1$.\\
 Particularly, if $\mu\{x\}=0$,
then $\grad f(x)= H(x)$. Moreover, $H$ is continuous outside the
support of $\mu$.
\end{proposition}

\begin{proof}
We prove only the first identity since the proof of the second one
is similar. For this, let $t_0\in [\,0,b\,)$ and
$\varepsilon\in(0,b-t_0]
$, then\\
\begin{align*}
&\frac{f(\gamma(t_0+\varepsilon))-f(\gamma(t_0))}{\varepsilon}
=\int_M\frac{d(\gamma(t_0+\varepsilon),p)-d(\gamma(t_0),p)}{\varepsilon}\mu(dp)\\
&=\int_{M\setminus\{\gamma(t_0)\}}\frac{d(\gamma(t_0+\varepsilon),p)-d(\gamma(t_0),p)}{\varepsilon}\mu(dp)+\mu\{\gamma(t_0)\}|\dot{\gamma}|
\end{align*}
By letting $\varepsilon\rightarrow0+$ and using the bounded
convergence to the above first integral we obtain that
\begin{align*}
\frac{d}{dt}f(\gamma(t))\big{|}_{t=t_0+}
&=\int_{M\setminus\{\gamma(t_0)\}}\frac{d}{dt}d(\gamma(t),p)\big{|}_{t=t_0}\mu(dp)+\mu\{\gamma(t_0)\}|\dot{\gamma}|\\
&=\int_{M\setminus\{\gamma(t_0)\}}\langle\,\dot{\gamma}(t_0),\,\grad
d_p(\gamma(t_0))\,\rangle\mu(dp)
+\mu\{\gamma(t_0)\}|\dot{\gamma}|\\
&=\langle\,\dot{\gamma}(t_0),\,H(\gamma(t_0))\,\rangle+\mu\{\gamma(t_0)\}|\dot{\gamma}|
\end{align*}

\end{proof}

Now we give a characterization of $\mathfrak{M}_{\mu}$ which is
proved in Khun \cite{Khun} for the case when $\mu$ is a finite set
of points in an Euclidean space.

\begin{theorem}\label{characterization of median}
$\mathfrak{M}_{\mu}=\big\{x\in
\bar{B}(a,\rho):|H(x)|\leq\mu{\{x\}}\big\}$
\end{theorem}

\begin{proof}
$(\,\subset\,)$ Let $x\in \mathfrak{M}_{\mu}$, if $H(x)=0$ there is
nothing to prove, so we assume that $H(x)\neq 0$. Consider a
geodesic in $\bar{B}(a,\rho)$:
\[\gamma(t)=\exp_x(-t\frac{H(x)}{|H(x)|}),\quad t\in [\,0,\,b\,]\]
according to the definition of $\mathfrak{M}_{\mu}$, $t=0$ is a
minimum point of $f\circ\gamma$, thus
\[\frac{d}{dt}f(\gamma(t))\big{|}_{t=0+}\geq 0\]
Moreover, by Proposition \ref{right derivative of f} we have
\[\frac{d}{dt}f(\gamma(t))\big{|}_{t=0+}=\langle\,-\frac{H(x)}{|H(x)|},\,H(x)\,\rangle+\mu\{x\}=-|H(x)|+\mu\{x\}\]
thus we have $|H(x)|\leq\mu{\{x\}}$.\\

$(\,\supset)$ Suppose that $|H(x)|\leq\mu{\{x\}}$, then for every
$y\in \bar{B}(a,\rho)$, we consider the geodesic
$\gamma:[0,1]\rightarrow \bar{B}(a,\rho)$ with $\gamma(0)=x$ and
$\gamma(1)=y$, then Proposition \ref{right derivative of f} and
Cauchy-Schwartz inequality give that
\begin{align*}
\frac{d}{dt}f(\gamma(t))\big{|}_{t=0+}&=\langle\,\dot{\gamma}(0),\,H(x)\,\rangle
+\mu\{x\}|\dot{\gamma}|\\
&\geq 
|\dot{\gamma}|(-|H(x)|+\mu\{x\})\geq 0
\end{align*}
By Lemma \ref{convexity of f}, $f\circ\gamma$ is convex then we have
\[f(y)-f(x)\geq\frac{d}{dt}f(\gamma(t))\big{|}_{t=0+}\geq 0\]
hence we get $x\in \mathfrak{M}_{\mu}$.\\
\end{proof}

In order to describe the location of the median of $\mu$, we need
the following geometric lemma which is also useful in the next
section.

\begin{lemma}\label{obtuse angle}
Let $\triangle ABC$ be a geodesic triangle in $\bar{B}(a,\rho)$ such
that $\angle A\geq \pi/2$, then $\angle B<\pi/2$ and $\angle
C<\pi/2$.
\end{lemma}

\begin{proof}

It suffices to prove that $\angle B<\pi/2$ . We show this only for
$\Delta>0$ since the cases when $\Delta\leq0$ are similar. Let
$d(B,C)=a$, $d(A,C)=b$ and $d(A,B)=c$. We consider a comparison
triangle $\triangle \bar{A}\bar{B}\bar{C}$ of $\triangle ABC$ in
$M^2_{\Delta}$, since $K\leq\Delta$ in $\bar{B}(a,\rho)$, hence by
Alexandrov's theorem we get $\angle A\leq\angle \bar{A}$ and $\angle
B\leq\angle \bar{B}$. The law of cosines in $M^2_{\Delta}$ together
with $\cos\angle \bar{A}\leq 0$ yields that
$\cos(\sqrt{\Delta}a)\leq \cos(\sqrt{\Delta}b)\cos(\sqrt{\Delta}c)$.
Using the law of cosines again we get
\begin{align*}
\cos\angle
\bar{B}&=\frac{\cos(\sqrt{\Delta}b)-\cos(\sqrt{\Delta}a)\cos(\sqrt{\Delta}c)}{\sin(\sqrt{\Delta}a)\sin(\sqrt{\Delta}c)}\\
&\geq\frac{\cos(\sqrt{\Delta}b)-\cos(\sqrt{\Delta}b)\cos^2(\sqrt{\Delta}c)}{\sin(\sqrt{\Delta}a)\sin(\sqrt{\Delta}c)}\\
&=\frac{\cos(\sqrt{\Delta}b)\sin(\sqrt{\Delta}c)}{\sin(\sqrt{\Delta}a)}>0
\end{align*}
thus $\angle B\leq\angle \bar{B}<\pi/2$.

\end{proof}

\begin{proposition}\label{location}
$\mathfrak{M}_{\mu}$ is contained in the smallest closed convex
subset of $B(a,\rho)$ containing the support of $\mu$.
\end{proposition}

\begin{proof}
Let $V$ be this set and by Theorem \ref{characterization of median}
it suffices to show that if $x\in\bar{B}(a,\rho)\setminus V$ then
$H(x)\neq 0$. In fact, let $y$ be a point in $V$ such that
$d(x,y)=\inf\{d(x,p):p\in V\}$, then by the convexity of $V$ we get
$\angle xyp\geq\pi/2$ for every $p\in V$ and hence Lemma \ref{obtuse
angle} yields that $\angle pxy<\pi/2$. This gives that
\begin{align*} \langle\,H(x),\exp^{-1}_xy\,\rangle
&=\int_{V}\frac{\langle\,-\exp^{-1}_{x}p,\,\exp^{-1}_xy\,\rangle}{d(x,p)}\mu(dp)\\
&=-d(x,y)\int_{V}\cos\angle pxy <0
\end{align*}
The proof is completed by observing that $\exp^{-1}_xy\neq0$.

\end{proof}

\section{uniqueness of riemannian median}

In the Euclidean case, if the sample points are not colinear, then
the geometric median is unique. Hence we get a natural condition of
$\mu$ to ensure the uniqueness of the median in Riemannian case:\\

\emph{$\ast$ \quad The support of $\mu$ is not totally contained in
any geodesic. This means that for every geodesic $\gamma$:
$[\,0,b\,]\rightarrow\bar{B}(a,\rho)$, we have
$\mu(\gamma[\,0,b\,])<1$.}\\

\begin{theorem}\label{uniqueness of meidan}
If condition $\ast$ holds, then $\mathfrak{M}_{\mu}$ has a single
element.
\end{theorem}

\begin{proof}
We will prove this by showing that $f$ is strictly convex, that is,
for every geodesic $\gamma$: $[\,0,b\,]\rightarrow\bar{B}(a,\rho)$,
the function $f\circ\gamma$ is strictly convex. In fact, without
loss of generality, we may assume that $\gamma(0)$ and $\gamma(b)$
are both in $\partial\bar{B}(a,\rho)$. By the first estimation in
Lemma \ref{hessian estimation}, for every
$p\in\bar{B}(a,\rho)\setminus\gamma[\,0,b\,]$ the function $t\mapsto
d(\gamma(t),p)$ is strictly convex and for $p\in\gamma[\,0,b\,]$ it
is trivially convex. Since the condition $\ast$ yields that
$\mu(\bar{B}(a,\rho)\setminus\gamma[\,0,b\,])>0$, thus by
integration we obtain the strict convexity of $f$ and this completes
the proof.
\end{proof}

In the above proof, we have seen that $f$ is strictly convex if
condition $\ast$ holds. However, things are better than this. In
fact, we can show that the condition $\ast$ implies the strong
convexity of $f$. In fact, the compactness of $\bar{B}(a,\rho)$ and
the condition $\ast$ can give a stronger result and the following
lemma clarifies this.

\begin{lemma}\label{constant of mu}

If condition $\ast$ holds, then there exist two constants
$\varepsilon_{\mu}\in(0,\rho)$ and $\eta_{\mu}\in (0,1]$ such that
for every geodesic $\gamma$: $[\,0,b\,]\rightarrow\bar{B}(a,\rho)$
we have
\[\mu(B(\gamma,\varepsilon_{\mu}))\leq 1-\eta_{\mu}\]
where for $\varepsilon>0$,
$B(\gamma,\varepsilon)=\{x\in\bar{B}(a,\rho): d(x,\gamma[\,0,b])<
\varepsilon\}$.

\end{lemma}

\begin{proof}
If this is not true, then for every $\varepsilon\in(0,\rho)$ and
$\eta\in(0,1]$, there exists a geodesic $\gamma$:
$[\,0,b\,]\rightarrow\bar{B}(a,\rho)$ such that
$\mu(B(\gamma,\varepsilon))> 1-\eta$. Then we obtain a sequence of
geodesics $(\gamma_n)_n$: $[\,0,b\,]\rightarrow\bar{B}(a,\rho)$
verifying $\mu(B(\gamma_n,1/n))> 1-1/n$ for sufficiently large $n$.
Since $\bar{B}(a,\rho)$ is compact, there exists a subsequence
$(\gamma_{n_k})_k$ and a geodesic $\gamma$:
$[\,0,b\,]\rightarrow\bar{B}(a,\rho)$ such that
$\gamma_{n_k}\rightarrow\gamma$ uniformly on $[\,0,b\,]$. Then for
every $j\geq1$, when $k$ is sufficiently large we have
$B(\gamma_{n_k},1/n_k)\subset B(\gamma,1/j)$, hence for these $k$ we
have $\mu(B(\gamma,1/j))\geq 1-1/n_k$, by letting
$k\rightarrow\infty$ we get $\mu(B(\gamma,1/j))=1$ and then by
letting $j\rightarrow\infty$ we get $\mu(\gamma[0,b\,])=1$ which
contradicts the condition $\ast$.
\end{proof}

\begin{lemma}\label{c one and c two}
Let $(\triangle A_iB_iC_i)_{i=1,2}$ be two geodesic triangles in
model space $M^2_{\kappa}$ such that $\angle A_1\leq \angle A_2$,
$\angle B_1\leq \pi/2$, $\bar{d}(A_1,C_1)=\bar{d}(A_2,C_2)$ and
$\bar{d}(A_1,B_1)=\bar{d}(A_2,B_2)$. Then $\angle C_2\leq\angle
C_1$.
\end{lemma}

\begin{proof}
We prove the lemma for the case when $\kappa>0$, the proof for the
cases when $\kappa\leq0$ is similar. Let $\bar{d}(B_1,C_1)=a_1$,
$\bar{d}(B_2,C_2)=a_2$, $\bar{d}(A_1,C_1)=\bar{d}(A_2,C_2)=b$ and
$\bar{d}(A_1,B_1)=\bar{d}(A_2,B_2)=c$. Since $\angle A_1\leq \angle
A_2$, hence $a_1\leq a_2$. By the law of cosines in $M^2_{\kappa}$,
we get
\begin{align*}
\frac{d}{da_2}\cos\angle
C_2&=\frac{d}{da_2}\frac{\cos(\sqrt{\kappa}c)-\cos(\sqrt{\kappa}a_2)\cos(\sqrt{\kappa}b)}{\sin(\sqrt{\kappa}b)\sin(\sqrt{\kappa}a_2)}\\
&=\frac{\sqrt{\kappa}(\cos(\sqrt{\kappa}b)-\cos(\sqrt{\kappa}a_2)\cos(\sqrt{\kappa}c))}{\sin(\sqrt{\kappa}b)\sin^2(\sqrt{\kappa}a_2)}\\
&\geq\frac{\sqrt{\kappa}(\cos(\sqrt{\kappa}b)-\cos(\sqrt{\kappa}a_1)\cos(\sqrt{\kappa}c))}{\sin(\sqrt{\kappa}b)\sin^2(\sqrt{\kappa}a_2)}\\
&=\frac{\sqrt{\kappa}\sin(\sqrt{\kappa}a_1)\sin(\sqrt{\kappa}c)}{\sin(\sqrt{\kappa}b)\sin^2(\sqrt{\kappa}a_2)}\cos
\angle B_1\geq 0
\end{align*}
so that $\cos\angle C_2$ is nondecreasing with respect to $a_2$ when
$a_2\geq a_1$, thus we get $\cos\angle C_2\geq\cos\angle C_1$, i.e.
$\angle C_2\leq\angle C_1$.
\end{proof}

\begin{lemma}\label{estimation of sinC}
Let $\triangle ABC$ be a geodesic triangle in $\bar{B}(a,\rho)$ such
that $\angle A= \pi/2$ then we have\\ \[\sin\angle C\geq
L(\rho,\delta)d(A,B)\] where the constant
\[
L(\rho,\delta)=
\begin{cases}
2\sqrt{\delta}\big/\big(\,\pi\sqrt{1-\cos^4(2\rho\sqrt{\delta})}\,\big)  & \text{if $\delta>0$;}\\
1\big/(2\sqrt{2}\rho)   & \text{if $\delta=0$;}\\
\sqrt{-\delta}\big/\sqrt{\cosh^4(2\rho\sqrt{-\delta})-1}  & \text{if
$\delta<0$;}
\end{cases}
\]
\end{lemma}

\begin{proof}
Let $(\triangle A_iB_iC_i)_{i=1,2}$ be two geodesic triangles in
$M^2_{\delta}$ such that $\triangle A_1B_1C_1$ is a comparison
triangle of $\triangle ABC$ and $\triangle A_2B_2C_2$ verifies that
$\bar{d}(A_2,C_2)=d(A,C)$, $\bar{d}(A_2,B_2)=d(A,B)$ and $\angle
A_2=\angle A$. Assume that $d(B,C)=a$, $d(A,C)=b$, $d(A,B)=c$ and
$\bar{d}(B_2,C_2)=a_2$. By Lemma \ref{obtuse angle} we get that
$\angle B<\pi/2$. Since $K\geq\delta$ in $\bar{B}(a,\rho)$, then by
Toponogov's theorem we get $\angle A_1\leq\angle A=\angle A_2$ and
$\angle B_1\leq\angle B<\pi/2$, hence by Lemma \ref{c one and c two}
we obtain that $\angle C_2\leq \angle C_1$. Observe that we have
also $\angle C_1\leq\angle C<\pi/2$, thus $\sin\angle  C\geq\sin
\angle C_2$. Now it suffices to estimate $\sin\angle  C_2$ according
to the sign of the lower sectional curvature bound $\delta$:\\

If $\delta>0$, using $\angle A_2=\pi/2$ and the laws of sines and
cosines in $M^2_{\delta}$ we get
\[\sin\angle  C_2=\frac{\sin(\sqrt{\delta}c)}{\sin(\sqrt{\delta}a_2)}
\quad\text{and}\quad\cos(\sqrt{\delta}a_2)=\cos(\sqrt{\delta}b)\cos(\sqrt{\delta}c)\]
hence we obtain that
\[\sin\angle  C_2=\frac{\sin(\sqrt{\delta}c)}{\sqrt{1-\cos^2(\sqrt{\delta}b)\cos^2(\sqrt{\delta}c)}}
\geq\frac{2\sqrt{\delta}c}{\pi\sqrt{1-\cos^4(2\rho\sqrt{\delta})}}\]
since $b,c\leq 2\rho$ and $\sin\theta\geq 2\theta/\pi$ for
$\theta\in[0,\pi/2]$.

Using the same method as above, we obtain that
\begin{align*}
&\text{If}\quad\delta=0,\quad\sin\angle
C_2=\frac{c}{\sqrt{b^2+c^2}}\geq\frac{c}{2\sqrt{2}\rho}\\
&\text{If}\quad\delta<0,\quad\sin\angle
C_2=\frac{\sinh(\sqrt{-\delta}c)}{\sqrt{\cosh^2(\sqrt{-\delta}b)\cosh^2(\sqrt{-\delta}c)-1}}
\geq\frac{\sqrt{-\delta}c}{\sqrt{\cosh^4(2\rho\sqrt{-\delta})-1}}
\end{align*}
since $\sinh\theta\geq\theta$ for $\theta\geq0$.
\end{proof}

We know that $f$ is convex, thus along every geodesic it has a
second derivative in the sense of distribution, the following
proposition gives its specific form as well as the Taylor's formula.

\begin{proposition}\label{Taylor}
Let $\gamma$: $[\,0,b\,]\rightarrow\bar{B}(a,\rho)$ be a geodesic
then for $t\in[\,0,b\,]$
\begin{align*}
f(\gamma(t))=f(\gamma(0))+\frac{d}{ds}f(\gamma(s))\big|_{s=0+}t+\int_{(0,t)}(t-s)\nu(ds)
\end{align*}
where $\nu$ is the second derivative of $f\circ\gamma$ on $(0,b)$ in
the sense of distribution, i.e. a bounded positive measure, which is
given by
\[\nu=\bigg(\int_{M\setminus\gamma[0,b\,]}\hess
d_p(\dot{\gamma},\dot{\gamma})\mu(dp)\bigg)\cdot\lambda|_{(0,b)}+2|\dot{\gamma}|\cdot(\mu\circ\gamma)|_{(0,b)}\]
with $\lambda|_{(0,b)}$ and $(\mu\circ\gamma)|_{(0,b)}$ denoting the
restrictions of Lebesgue measure and the measure $\mu\circ\gamma$ on
$(0,b)$ respectively.


\end{proposition}

\begin{proof}
Firstly, observe that since $\gamma$ is a homeomorphism of $(0,b)$
onto its image, $\mu\circ\gamma$ is a well defined measure on
$(0,b)$. By Proposition \ref{right derivative of f} we get
\begin{align*}
&\int_{M\setminus\gamma[0,b]}d(\gamma(t),p)\mu(dp)\\
=&\int_{M\setminus\gamma[0,b]}\bigg(d(\gamma(0),p)+\frac{d}{ds}d(\gamma(s),p)\big|_{s=0}t
+\int_0^t(t-s)\frac{d^2}{ds^2}d(\gamma(s),p)ds\bigg)\mu(dp)\\
=&\int_{M\setminus\gamma[0,b]}d(\gamma(0),p)\mu(dp)
+\langle\,\dot{\gamma}(0),\int_{M\setminus\gamma[0,b]}\frac{-\exp^{-1}_{\gamma(0)}p}
{d(\gamma(0),p)}\mu(dp)\,\rangle t\\
&+\int_0^t(t-s)ds\int_{M\setminus\gamma[0,b]}\hess
d_p(\dot{\gamma}(s),\dot{\gamma}(s))\mu(dp)\\
=& \,f(\gamma(0))-\int_{\gamma(0,b]}d(\gamma(0),p)\mu(dp)\\
&+\frac{d}{ds}f(\gamma(s))\big|_{s=0+}t
-\langle\,\dot{\gamma}(0),\int_{\gamma(0,b]}\frac{-\exp^{-1}_{\gamma(0)}p}
{d(\gamma(0),p)}\mu(dp)\,\rangle t-\mu\{\gamma(0)\}|\dot{\gamma}|t\\
&+\int_0^t(t-s)ds\int_{M\setminus\gamma[0,b]}\hess
d_p(\dot{\gamma}(s),\dot{\gamma}(s))\mu(dp)
\end{align*}
Since
\[f(\gamma(t))=\int_{M\setminus\gamma[0,b]}d(\gamma(t),p)\mu(dp)+\int_{\gamma[0,b]}d(\gamma(t),p)\mu(dp))\]
we obtain
\begin{align*}
&f(\gamma(t))-f(\gamma(0))-\frac{d}{ds}f(\gamma(s))\big|_{s=0+}t-\int_0^t(t-s)ds\int_{M\setminus\gamma[0,b]}\hess
d_p(\dot{\gamma}(s),\dot{\gamma}(s))\mu(dp)\\
&=-\int_{\gamma(0,b]}d(\gamma(0),p)\mu(dp)+\mu(\gamma(0,b])|\dot{\gamma}|t-\mu\{\gamma(0)\}|\dot{\gamma}|t+\int_{\gamma[0,b]}d(\gamma(t),p)\mu(dp)\\
&=-\int_{\gamma(0,b]}d(\gamma(0),p)\mu(dp)+\int_{\gamma(0,b]}d(\gamma(0),\gamma(t))\mu(dp)+\int_{\gamma(0,b]}d(\gamma(t),p)\mu(dp)\\
&=\bigg(-\int_{\gamma(0,t)}d(\gamma(0),p)\mu(dp)+\int_{\gamma(0,t)}d(\gamma(0),\gamma(t))\mu(dp)+\int_{\gamma(0,t)}d(\gamma(t),p)\mu(dp)\bigg)\\
&+\bigg(-\int_{\gamma[t,b]}d(\gamma(0),p)\mu(dp)+\int_{\gamma[t,b]}d(\gamma(0),\gamma(t))\mu(dp)+\int_{\gamma[t,b]}d(\gamma(t),p)\mu(dp)\bigg)\\
&=2\int_{\gamma(0,t)}d(\gamma(t),p)\mu(dp)+0=2|\dot{\gamma}|\int_{(0,t)}(t-s)(\mu\circ\gamma)(ds)
\end{align*}
hence the desired formula holds. To show that $\nu$ is the second
derivative of $f\circ\gamma$ on $(0,b)$ in the sense of
distribution, let $\varphi\in C^{\infty}_c(\,0,b)$ and by Fubini's
theorem and integration by parts we get
\begin{align*}
&\int_{(0,b)}f(\gamma(t))\varphi^{\prime\prime}(t)dt\\
=&\,f(\gamma(0))\int_{(0,b)}\varphi^{\prime\prime}(t)dt
+\frac{d}{ds}f(\gamma(s))\big|_{s=0+}\int_{(0,b)}t\varphi^{\prime\prime}(t)dt
+\int_{(0,b)}\varphi^{\prime\prime}(t)dt\int_{(0,t)}(t-s)\nu(ds)\\
=&\int_{(0,b)}\nu(ds)\int_{(s,b)}(t-s)\varphi^{\prime\prime}(t)dt=\int_{(0,b)}\varphi(s)\nu(ds)
\end{align*}
the proof is complete.
\end{proof}

Now we are ready to show the strong convexity of $f$ under condition
$\ast$ which will be useful for our error estimates. Certainly, this
also yields the uniqueness of the median.

\begin{theorem}\label{strong convex}
If condition $\ast$ holds, then for every geodesic $\gamma$:
$[\,0,b\,]\rightarrow\bar{B}(a,\rho)$, we have the following
inequality
\[f(\gamma(t))\geq
f(\gamma(0))+\frac{d}{ds}f(\gamma(s))\big|_{s=0+}t+\tau|\dot{\gamma}|^2
t^2,\quad t\in [\,0,b\,]\] where the constant
$\tau=(1/2)\,\varepsilon^2_{\mu}\,\eta_{\mu}\,D(\rho,\Delta)\,L(\rho,\delta)^2$.\\
Hence $f$ is $\tau$-convex. That is, for every unit velocity
geodesic $\gamma$: $[\,0,b\,]\rightarrow\bar{B}(a,\rho)$, the
function $t\mapsto f(\gamma(t))-\tau t^2$ is convex.
\end{theorem}

\begin{proof}
Without loss of generality we may assume that $\gamma(0)$ and
$\gamma(b)$ are both in $\partial\bar{B}(a,\rho)$. Then by the first
estimation in Lemma \ref{hessian estimation} we obtain that for
every $s\in [\,0,b\,]$,
\begin{align*}
&\int_{M\setminus\gamma[0,b]}\hess
d_p(\dot{\gamma}(s),\dot{\gamma}(s))\mu(dp)\geq\int_{M\setminus\gamma[0,b]}D(\rho,\Delta)|\dot{\gamma}^{\nor}(s)|^2\mu(dp)\\
=&\,D(\rho,\Delta)|\dot{\gamma}|^2\int_{M\setminus
\gamma[0,b]}\sin^2\angle(\dot{\gamma}(s),\exp^{-1}_{\gamma(s)}p)\mu(dp)
\end{align*}
Then for every $p\in M\setminus\gamma[\,0,b\,]$, let $q=q(p)$ be the
orthogonal projection of $p$ onto $\gamma$. If $\gamma(s)\neq q$,
the triangle $\triangle pq\gamma(s)$ is a right triangle with
$\angle pq\gamma(s)=\pi/2$. Hence Lemma \ref{estimation of sinC}
yields that
\[\sin\angle(\dot{\gamma}(s)),\exp^{-1}_{\gamma(s)}p)\geq L(\rho,\delta)d(p,q)\]
thus by Lemma \ref{constant of mu} we obtain
\begin{align*}
&\int_{M\setminus\gamma[0,b\,]}\hess
d_p(\dot{\gamma}(s),\dot{\gamma}(s))\mu(dp)\\
\geq\,&D(\rho,\Delta)|\dot{\gamma}|^2L(\rho,\delta)^2\int_{M\setminus\gamma[0,b\,]}d^2(p,q)\mu(dp)\\
\geq\,&D(\rho,\Delta)|\dot{\gamma}|^2L(\rho,\delta)^2\int_{M\setminus
B(\gamma,\varepsilon_{\mu})}d^2(p,q)\mu(dp)\\
\geq\,&\varepsilon^2_{\mu}\,\eta_{\mu}\,
D(\rho,\Delta)\,L(\rho,\delta)^2|\dot{\gamma}|^2=2\tau|\dot{\gamma}|^2
\end{align*}
Then by Proposition \ref{Taylor} we get that for $t\in [\,0,b\,]$,
\begin{align*}
f(\gamma(t))&=f(\gamma(0))+\frac{d}{ds}f(\gamma(s))\big|_{s=0+}t+\int_{(0,t)}(t-s)\nu(ds)\\
&\geq f(\gamma(0))+\frac{d}{ds}f(\gamma(s))\big|_{s=0+}t+2\,\tau|\dot{\gamma}|^2\int_{(0,t)}(t-s)ds\\
&\geq
f(\gamma(0))+\frac{d}{ds}f(\gamma(s))\big|_{s=0+}t+\tau|\dot{\gamma}|^2
t^2
\end{align*}
hence the desired inequality holds. To show the $\tau$-convexity of
$f$, let $\gamma$: $[\,0,b\,]\rightarrow\bar{B}(a,\rho)$ be a unit
velocity geodesic, then the above inequality shows that the function
$f(\gamma(t))-\tau t^2$ has an affine support on every $t\in
[\,0,b\,]$, so that it is convex.
\end{proof}

\section{a subgradient algorithm}
To begin with, we recall the notion of subgradient for a convex
function on a Riemannian manifold. For our purpose, it suffices to
consider this notion in a convex subset of the manifold.

\begin{definition}
Let $U$ be a convex subset of $M$ and $h$ be a convex function
defined on $U$. For every $x\in U$,  a vector $v\in T_xM$ is called
a subgradient of $h$ at $x$ if for every geodesic $\gamma$:
$[\,0,b\,]\rightarrow U$ with $\gamma(0)=x$, we have for every $t\in
[\,0,b\,]$,
\[h(\gamma(t))\geq h(x)+\,\langle\,\dot{\gamma}(0),\,v\,\rangle\,t\]

\end{definition}

Our idea to approximate the Riemannian median of $\mu$ by the
subgradient  method stems from the following simple observation.
\begin{lemma}\label{h is subgradient of f}
For every $x\in\bar{B}(a,\rho)$, $H(x)$ is a subgradient of $f$ at
$x$.
\end{lemma}

\begin{proof}
Let $\gamma$: $[\,0,1]\rightarrow\bar{B}(a,\rho)$ be a geodesic such
that $\gamma(0)=x$, then by Proposition \ref{right derivative of f}
together with the convexity of $f$ we get for every $t\in [\,0,1\,]
$,
\begin{align*} f(\gamma(t))&\geq
f(\gamma(0))+\frac{d}{ds}f(\gamma(s))\big|_{s=0+}t\\
&=f(x)+(\,\langle\,\dot{\gamma}(0),\,H(x)\,\rangle+\mu\{x\}|\dot{\gamma}|\,)\,t\\
&\geq f(x)+\,\langle\,\dot{\gamma}(0),\,H(x)\,\rangle \, t
\end{align*}
thus the assertion holds.
\end{proof}

In order to introduce our subgradient algorithm we need the
following notations.

\begin{notation}

Let $x\in \bar{B}(a,\rho)$ with $H(x)\neq 0$, then we write
\begin{align*}
\gamma_{x}(t)&=\exp_x(-t\frac{H(x)}{|H(x)|})\,,\quad t\geq 0\\
r_x&=\sup\{t\in[\,0,2\rho\,]:\gamma_{x}(t)\in \bar{B}(a,\rho)\}
\end{align*}

\end{notation}

We can now describe our subgradient algorithm to approximate the
Riemannian median of the probability measure $\mu$.

\begin{algorithm}\label{algo1}Subgradient algorithm for Riemannian median:\\
\textbf{Step 1:}\\
Choose a point $x_0\in \bar{B}(a,\rho)$ and let $k=0$.\\
\textbf{Step 2:}\\
If $|H(x_k)|\leq \mu\{x_k\}$, then stop and let $m=x_k$.\\
If not, then go to step 3.\\
\textbf{Step 3:}\\
Choose a stepsize $t_k\in(\,0,\,r_{x_k}]$ and let
$x_{k+1}=\gamma_{x_k}(t_k)$, then come back to step 2 with $k=k+1$.

\end{algorithm}

Observe that according to the definition of $r_x$, the sequence
$(x_k)_k$ is contained in the ball $\bar{B}(a,\rho)$. Now we turn to
the convergence proof of the above algorithm under some conditions
of the stepsize. It is well known that the following type of
inequalities are of fundamental importance to conclude the
convergence of the subgradient algorithms in Euclidean spaces:
\[||x_{k+1}-y||^2\leq ||x_k-y||^2+A_1t^2_k+A_2\frac{2t_k}{||v_k||}(f(y)-f(x_k))\]
see for example, Correa and Lemar\'{e}chal \cite{Correa}, Nedic and
Bertsekas \cite{Nedic2}. For a positively curved Riemannian
manifold, Ferreira and Oliveira \cite{Ferreira} obtained a
generalization of the above inequality by using Toponogov's
comparison theorem. But their method is not applicable in our case
since we do not suppose that $\delta\geq0$. However, we can still
obtain a similar result using a different method.

\begin{lemma}\label{fundamental lemma}
If $H(x_k)\neq0$, then for every point $y\in \bar{B}(a,\rho)$ we
have
\[d^2(x_{k+1},y)\leq d^2(x_k,y)+C(\rho,\delta)t^2_k+\frac{2t_k}{|H(x_k)|}(f(y)-f(x_k))\]
Particularly,
\[d^2(x_{k+1},\mathfrak{M}_{\mu})\leq d^2(x_k,\mathfrak{M}_{\mu})+C(\rho,\delta)t^2_k+2t_k(f_*-f(x_k))\]
\end{lemma}

\begin{proof}
By Taylor's formula and the second estimation in Lemma \ref{hessian
estimation}, there exists $\xi\in(0,r_x)$ such that
\begin{align*}
\frac{1}{2}d^2(x_{k+1},y)
&=\frac{1}{2}d^2(\gamma_{x_k}(t_k),y)\\
&=\frac{1}{2}d^2(\gamma_{x_k}(0),y)+\frac{d}{dt}\bigg[\frac{1}{2}d^2(\gamma_{x_k}(t),y)\bigg]_{t=0}t_k+\frac{1}{2}\frac{d^2}{dt^2}\bigg[\frac{1}{2}d^2(\gamma_{x_k}(t),y)\bigg]_{t=\xi}t^2_k\\
&=\frac{1}{2}d^2(x_k,y)+\langle\,\dot{\gamma}_{x_k}(0),\,\grad
\frac{1}{2}d^2_y(x_k)\,\rangle\, t_k+\frac{1}{2}\hess
\frac{1}{2}d^2_y(\dot{\gamma}(\xi),\dot{\gamma}(\xi))
\,t^2_k\\
&\leq\frac{1}{2}d^2(x_k,y)+\frac{\langle\,H(x_k),\,\exp^{-1}_{x_k}
y\,\rangle}{|H(x_k)|}\,t_k+\frac{C(\rho,\delta)}{2}\,t^2_k
\end{align*}
By Lemma \ref{h is subgradient of f}, $H(x_k)$ is a subgradient of
$f$ at point $x_k$ and hence
\[\langle\,H(x_k),\,\exp^{-1}_{x_k}
y\,\rangle\leq f(y)-f(x_k)\] thus we get
\[\frac{1}{2}d^2(x_{k+1},y)\leq\frac{1}{2}d^2(x_k,y)+\frac{t_k}{|H(x_k)|}(f(y)-f(x_k))+\frac{C(\rho,\delta)}{2} t^2_k\]
then the first inequality holds. The second one follows from
$f_*\leq f(x_k)$ and $|H(x_k)|\leq1$.
\end{proof}





As in the Euclidean case, once the fundamental inequality is
established, the convergence of the subgradient algorithm is soon
achieved, see for example Correa and Lemar\'{e}chal \cite{Correa},
Nedic and Bertsekas \cite{Nedic2}. Since our fundamental inequality
is very similar to the one in Euclidean case, the proof of
convergence is also very similar.

\begin{theorem}\label{convergence1}
If the stepsize $(t_k)_k$ verifies
 \[\lim_{k\rightarrow\infty}t_k=0\quad\text{and}\quad\sum_{k=0}^{\infty}t_k=+\infty\]
then we have
\[\lim_{k\rightarrow\infty}d(x_k,\mathfrak{M}_{\mu})=0\quad\text{and}\quad\lim_{k\rightarrow\infty}f(x_k)=f_*\]
Moreover, if the stepsize $(t_k)_k$ verifies also that
\[\sum_{k=0}^{\infty}t^2_k<+\infty\]
then there exists some $m\in \mathfrak{M}_{\mu}$ such that
$x_k\rightarrow m$, when $k\rightarrow\infty$.
\end{theorem}

\begin{proof}
If $|H(x_k)|\leq \mu\{x_k\}$ for some $k\geq0$, then Theorem
\ref{characterization of median} yields that $x_k\in
\mathfrak{M}_{\mu}$ and the the desired result is trivially true.
Hence in the following we assume that $|H(x_k)|> \mu\{x_k\}$ for
every $k$, this means that none of $x_k$ is in the median of $\mu$.
We will show firstly that
\[\liminf_{k\rightarrow\infty} f(x_k)=f_*\] If this is not true,
then there exist $N_1\in \mathbf{N}$ and $\eta>0$ such that for
every $k\geq N_1$ we have $f_*-f(x_k)\leq -\eta$.  Since
$|H(x_k)|>0$, then Lemma \ref{fundamental lemma} gives that
\begin{align*}
d^2(x_{k+1},\mathfrak{M}_{\mu}) \leq
d^2(x_k,\mathfrak{M}_{\mu})+t_k(C(\rho,\delta)t_k-2\eta)
\end{align*}
Since $\lim_{k\rightarrow\infty}t_k=0$, we can suppose that
$C(\rho,\delta)t_k<\eta$ for every $k\geq N_1$ and then
\[d^2(x_{k+1},\mathfrak{M}_{\mu})\leq d^2(x_k,\mathfrak{M}_{\mu})-\eta t_k\]
by summing the above inequalities we get
\[\eta\sum_{i=N_1}^{k}t_i\leq d^2(x_{N_1},\mathfrak{M}_{\mu})-d^2(x_{k+1},\mathfrak{M}_{\mu})\leq d^2(x_{N_1},\mathfrak{M}_{\mu})\]
which contradicts with $\sum_{k=0}^{\infty}t_k=+\infty$, this proves
the assertion.

Now for fixed $\varepsilon>0$, there exists $N_2\in\mathbf{N}$ such
that $C(\rho,\delta)t_k<2\varepsilon$ for every $k\geq N_2$. We consider the following two cases:\\
If $f(x_k)>f_*+\varepsilon$, then by Lemma \ref{fundamental lemma}
we obtain that
\begin{align*}
d^2(x_{k+1},\mathfrak{M}_{\mu})
&\leq d^2(x_{k},\mathfrak{M}_{\mu})+C(\rho,\delta)t^2_k+2t_k(f_*-f(x_k))\\
& <
d^2(x_{k},\mathfrak{M}_{\mu})+(C(\rho,\delta)t_k-2\varepsilon)t_k
<d^2(x_{k},\mathfrak{M}_{\mu})
\end{align*}
If $f(x_k)\leq f_*+\varepsilon$ then $x_k\in
L_{\varepsilon}=\{x\in\bar{B}(a,\rho):f(x)\leq f_*+\varepsilon\}$
and if we write
$l_{\varepsilon}=\sup\{d(\,y,\mathfrak{M}_{\mu}):y\in
L_{\varepsilon}\}$, hence in this case we have
\[d(x_{k+1},\mathfrak{M}_{\mu})\leq
d(x_{k+1},x_k)+d(x_k,\mathfrak{M}_{\mu})\leq t_k +l_{\varepsilon}\]
In conclusion, we always have that for $k\geq N_2$,
\[d(x_{k+1},\mathfrak{M}_{\mu})\leq\max\{d(x_{k},\mathfrak{M}_{\mu}),t_k+l_{\varepsilon}\}\]
By induction we get for every $n\geq k$,
\begin{align*}
d(x_{n+1},\mathfrak{M}_{\mu})
&\leq\max\{d(x_{k},\mathfrak{M}_{\mu}),\max\{t_k,t_{k+1},\dots,t_n\}+l_{\varepsilon}\}\\
&\leq\max\{d(x_{k},\mathfrak{M}_{\mu}),\sup\{t_i:i\geq
k\}+l_{\varepsilon}\}
\end{align*}
thus we get
\[\limsup_{n\rightarrow\infty}d(x_{n},\mathfrak{M}_{\mu})\leq\max\{d(x_{k},\mathfrak{M}_{\mu}),\sup\{t_i:i\geq
k\}+l_{\varepsilon}\}\] $\liminf_{k\rightarrow\infty} f(x_k)=f_*$
yields that $\liminf_{k\rightarrow\infty}
d(x_{k},\mathfrak{M}_{\mu})=0$, by taking the inferior limit  on the
right hand side we obtain
\[\limsup_{n\rightarrow\infty}d(x_{n},\mathfrak{M}_{\mu})\leq l_{\varepsilon}\]
Now we show that $l_{\varepsilon}\rightarrow0$ when
$\varepsilon\rightarrow0$. By monotonicity of $l_{\varepsilon}$, it
suffices to show this along some sequence. In fact, observe that
$L_{\varepsilon}$ is compact, thus for every $\varepsilon>0$, there
exists $y_{\varepsilon}\in L_{\varepsilon}$ such that
$l_{\varepsilon}=d(y_{\varepsilon},\mathfrak{M}_{\mu})$. Since
$\bar{B}(a,\rho)$ is compact, there exist a sequence
$\varepsilon_{k}\rightarrow 0$ and $y\in\bar{B}(a,\rho)$ such that
$y_{\varepsilon_k}\rightarrow y$. Since $f(y_{\varepsilon_k})\leq
f_*+\varepsilon_k$, we have $f(y)\leq f_*$ and hence $y\in
\mathfrak{M}_{\mu}$. Consequently $l_{\varepsilon_k}\rightarrow0$.
Thus we get $d(x_{k},\mathfrak{M}_{\mu})\rightarrow0$ and this
yields $f(x_k)\rightarrow
f_*$.\\
If $\sum_{k=0}^{\infty}t^2_k<+\infty$, the compactness of
$\bar{B}(a,\rho)$ and $f(x_k)\rightarrow f_*$ imply that the
sequence $(x_k)_k$ has some cluster point $m\in \mathfrak{M}_{\mu}$,
hence Lemma \ref{fundamental lemma} yields
\[d^2(x_{k+1},m)\leq d^2(x_k,m)+C(\rho,\delta)t^2_k\]
for every $n\geq k$, by summing the above inequalities we get
\[d^2(x_{n+1},m)\leq d^2(x_k,m)+C(\rho,\delta)\sum_{i=k}^{n}t^2_i\]\
and by letting $n\rightarrow\infty$ we obtain
\[\limsup_{n\rightarrow\infty}d^2(x_{n},m)\leq d^2(x_k,m)+C(\rho,\delta)\sum_{i=k}^{\infty}t^2_i\]\
the proof will be completed by observing that the right hand side of
the above inequality posseses a subsequence that converges to $0$.

\end{proof}

Now we have to consider the choice of stepsize. In fact, we can
choose $(t_k)_k$ that verifies the conditions of the preceding
theorem thus yields the desired convergence of our algorithm and
this is justified by the following lemma.

\begin{lemma}
$\inf\{r_x:x\in \bar{B}(a,\rho)\}>0$
\end{lemma}

\begin{proof}
Since the support of $\mu$ is contained in $B(a,\rho)$, if
$x\in\partial\bar{B}(a,\rho)$, $H(x)$ is transverse to
$\partial\bar{B}(a,\rho)$ and hence $r_x>0$ for every
$x\in\bar{B}(a,\rho)$. Moreover, there exists $\varepsilon>0$ such
that $\supp(\mu)\subset B(a,\rho-\varepsilon)$, then for $x\in
B(a,\rho-\varepsilon)$, we have $r_x\geq \rho-d(x,a)> \varepsilon$.
On the other hand, $H$ is continuous on $\bar{B}(a,\rho)\setminus
B(a,\rho-\varepsilon)$ which is compact, thus there exists a point
$x_0\in \bar{B}(a,\rho)\setminus B(a,\rho-\varepsilon)$ such that
$\text{inf}\{r_x: x\in \bar{B}(a,\rho)\setminus
B(a,\rho-\varepsilon)\}=r_{x_0}>0$. Hence we get $\inf\{r_x:x\in
\bar{B}(a,\rho)\}\geq \min\{\varepsilon,r_{x_0}\}>0$.
\end{proof}

According to this lemma, we may take for example the stepsize
$t_k=r_{x_k}/(k+1)$ for every $k\geq0$, then by Theorem
\ref{convergence1} the convergence holds. But the drawback is that
we do not know $r_{x_k}$. However, with further analysis we can give
an explicit uniform lower bound of them. In the following, let
\[\sigma=\sup\{d(p,a):p\in \supp\mu\}\] and observe that
$\supp(\mu)\subset B(a,\rho)$ yields $\sigma<\rho$.

\begin{lemma}
For every $x\in \bar{B}(a,\rho)\setminus\bar{B}(a,\sigma)$ we have

\[r_x \geq \frac{2d(x,a)S_{\Delta}(d(x,a)-\sigma)}{C(\rho,\delta)S_{\Delta}(d(x,a)+\sigma)}\]
\end{lemma}

\begin{proof}
Since the diameter of $\bar{B}(a,\rho)$ is $2\rho$, then
$\gamma_x(r_x)\in
\partial \bar{B}(a,\rho)$. By Taylor's formula and the second estimation in Lemma \ref{hessian
estimation}, there exists $\xi\in(0,r_x)$ such that
\begin{align*}
\frac{1}{2}\rho^2
&=\frac{1}{2}d^2(\gamma_x(r_x),a)\\
&=\frac{1}{2}d^2(x,a)+\frac{d}{dt}\bigg[\frac{1}{2}d^2(\gamma_x(t),a)\bigg]_{t=0}r_x+\frac{1}{2}\frac{d^2}{dt^2}\bigg[\frac{1}{2}d^2(\gamma_x(t),a)\bigg]_{t=\xi}r^2_x\\
&=\frac{1}{2}d^2(x,a)+\langle\,\dot{\gamma}_x(0),\,\grad
\frac{1}{2}d^2_a(x)\,\rangle\, r_x+\frac{1}{2}\hess
\frac{1}{2}d^2_a(\dot{\gamma}(\xi),\dot{\gamma}(\xi))
\,r^2_x\\
&\leq\frac{1}{2}d^2(x,a)+\frac{\langle\,H(x),\,\exp^{-1}_x
a\,\rangle}{|H(x)|}\,r_x+\frac{C(\rho,\delta)}{2}\,r^2_x
\end{align*}
Gauss Lemma yields that $\langle\,\exp^{-1}_x p,\,\exp^{-1}_x
a\,\rangle>0 $ for $p\in \supp\mu$, hence
\[\langle\,H(x),\,\exp^{-1}_x
a\,\rangle=-\int_{\supp\mu}\frac{\langle\,\exp^{-1}_x
p,\,\exp^{-1}_x a\,\rangle}{d(x,p)}\mu(dp)< 0\] Combine this with
$d(x,a)\leq\rho$, $C(\rho,\delta)>0$ and $r_x>0$, we obtain that
\begin{align*}
r_x&\geq\frac{1}{C(\rho,\delta)}\bigg\{-\frac{\langle\,H(x),\,\exp^{-1}_x
a\,\rangle}{|H(x)|}+\sqrt{\frac{\langle\,H(x),\,\exp^{-1}_x
a\,\rangle^2}{|H(x)|^2}+C(\rho,\delta)(\rho^2-d^2(x,a))}\bigg\}\\
&\geq\frac{1}{C(\rho,\delta)}\bigg\{-\frac{\langle\,H(x),\,\exp^{-1}_x
a\,\rangle}{|H(x)|}+\frac{|\langle\,H(x),\,\exp^{-1}_x
a\,\rangle|}{|H(x)|}\bigg\}\\
&=\frac{-2}{C(\rho,\delta)}\frac{\langle\,H(x),\,\exp^{-1}_x
a\,\rangle}{|H(x)|}\geq
\frac{-2}{C(\rho,\delta)}\langle\,H(x),\,\exp^{-1}_x a\,\rangle\\
&=\frac{2}{C(\rho,\delta)}\int_{\supp\mu}\frac{\langle\,\exp^{-1}_x
p,\,\exp^{-1}_x a\,\rangle}{d(x,p)}\mu(dp)\\
&=\frac{2d(x,a)}{C(\rho,\delta)}\int_{\supp\mu}\cos\angle\, pxa\,\,
\mu(dp)
\end{align*}
For every $p\in\supp(\mu)$ we consider the triangle $\triangle pxa$
and its comparison triangle $\triangle \bar{p}\bar{x}\bar{a}$ in
$M_{\Delta}^2$. Since $K\leq\Delta$ in $\bar{B}(a,\rho)$, then
Alexandrov's Theorem implies that $\angle\, pxa\leq \angle\,
\bar{p}\bar{x}\bar{a}$ and hence $\cos\angle\, pxa\geq\cos\angle\,
\bar{p}\bar{x}\bar{a}$. Now it suffices to estimate the lower bound
of $\cos\angle\, \bar{p}\bar{x}\bar{a}$ according to the sign of
$\Delta$. In fact, we have for every $\Delta\in\mathbf{R}$,
\[\cos\angle\,\bar{p}\bar{x}\bar{a}\geq\frac{S_{\Delta}(d(x,a)-\sigma)}{S_{\Delta}(d(x,a)+\sigma)}\]
We show this for the case when $\Delta>0$, the proof for cases when
$\Delta\leq0$ is similar. For this, let us observe that
$d(a,p)\leq\sigma$ and that $d(x,a)-\sigma\leq d(x,p)\leq
d(x,a)+\sigma$, then we obtain
\begin{align*}
\cos\angle\,\bar{p}\bar{x}\bar{a}
&=\frac{\cos(\sqrt{\Delta}d(a,p))-\cos(\sqrt{\Delta}d(x,a))\cos(\sqrt{\Delta}d(x,p))}{\sin(\sqrt{\Delta}d(x,a))\sin(\sqrt{\Delta}d(x,p))}\\
&\geq\frac{\cos(\sqrt{\Delta}\sigma)-\cos(\sqrt{\Delta}d(x,a))\cos(\sqrt{\Delta}(d(x,a)-\sigma))}{\sin(\sqrt{\Delta}d(x,a))\sin(\sqrt{\Delta}(d(x,a)+\sigma))}\\
&=\frac{\sin(\sqrt{\Delta}(d(x,a)-\sigma))}{\sin(\sqrt{\Delta}(d(x,a)+\sigma))}
=\frac{S_{\Delta}(d(x,a)-\sigma)}{S_{\Delta}(d(x,a)+\sigma)}
\end{align*}
hence the claimed inequality is proved and consequently,
\[r_x
\geq\frac{2d(x,a)}{C(\rho,\delta)}\int_{\supp\mu}\frac{S_{\Delta}(d(x,a)-\sigma)}{S_{\Delta}(d(x,a)+\sigma)}\mu(dp)
=\frac{2d(x,a)S_{\Delta}(d(x,a)-\sigma)}{C(\rho,\delta)S_{\Delta}(d(x,a)+\sigma)}\]
\end{proof}

Now we can give the desired lower bound.

\begin{lemma}\label{estimation of rx}

For every $x\in \bar{B}(a,\rho)$ we have
\[r_x\geq \frac{\rho-\sigma}{C(\rho,\delta)\cosh(2\rho\sqrt{|\Delta|})+1}\]

\end{lemma}

\begin{proof}
Firstly, assume that $x\in
\bar{B}(a,\rho)\setminus\bar{B}(a,\sigma)$ then by the preceding
lemma, if $\Delta>0$ we have
\[r_x\geq\frac{2d(x,a)}{C(\rho,\delta)}\frac{\sin(\sqrt{\Delta}(d(x,a)-\sigma))}{\sin(\sqrt{\Delta}(d(x,a)+\sigma))}\]
observe that
$0<\sqrt{\Delta}(d(x,a)+\sigma))<2\rho\sqrt{\Delta}\leq\pi/2$ and
that for $0<u\leq v\leq \pi/2$, we have $(\sin u/\sin v)\geq (u/v)$,
then we obtain
\[r_x
\geq\frac{2d(x,a)}{C(\rho,\delta)}\frac{d(x,a)-\sigma}{d(x,a)+\sigma}
\geq\frac{2d(x,a)}{C(\rho,\delta)}\frac{d(x,a)-\sigma}{2d(x,a)}
=\frac{d(x,a)-\sigma}{C(\rho,\delta)}\] On the other hand, we always
have $r_x\geq\rho-d(x,a)$ and hence
\[r_x\geq \max\{\rho-d(x,a),\frac{d(x,a)-\sigma}{C(\rho,\delta)}\}\]
Observe that
\[\min\bigg\{\max\{\rho-d(x,a),\frac{d(x,a)-\sigma}{C(\rho,\delta)}\}:\sigma<d(x,a)\leq\rho\bigg\}
=\frac{\rho-\sigma}{C(\rho,\delta)+1}\] hence we get
\[r_x\geq\frac{\rho-\sigma}{C(\rho,\delta)+1}\]
If $\Delta=0$, then we have
\[r_x\geq\frac{2d(x,a)}{C(\rho,\delta)}\frac{d(x,a)-\sigma}{d(x,a)+\sigma}\]
the same proof as above yields the same result.\\
If $\Delta<0$ then we have
\[r_x\geq\frac{2d(x,a)}{C(\rho,\delta)}\frac{\sinh(\sqrt{-\Delta}(d(x,a)-\sigma))}{\sinh(\sqrt{-\Delta}(d(x,a)+\sigma))}\]
Observe that for $u>0$, $\sinh u$ et $\cosh u$ are nondecreasing and
that $u\leq\sinh u\leq u\cosh u$, then we obtain
\[r_x\geq\frac{2d(x,a)}{C(\rho,\delta)}\frac{\sqrt{-\Delta}(d(x,a)-\sigma)}{\sinh(\sqrt{-\Delta}2d(x,a))}
\geq\frac{d(x,a)-\sigma}{C(\rho,\delta)\cosh(2\rho\sqrt{-\Delta})}\]
the same method as in the case $\Delta>0$ yields
\[r_x\geq\frac{\rho-\sigma}{C(\rho,\delta)\cosh(2\rho\sqrt{-\Delta})+1}\]
In conclusion, since $\cosh(2\rho\sqrt{|\Delta|})\geq 1$, we have
always
\[r_x\geq
\frac{\rho-\sigma}{C(\rho,\delta)\cosh(2\rho\sqrt{|\Delta|})+1}\]
for every $x\in\bar{B}(a,\rho)\setminus\bar{B}(a,\sigma)$. Moreover,
for those $x\in\bar{B}(a,\sigma)$,
\[r_x\geq\rho-\sigma>
\frac{\rho-\sigma}{C(\rho,\delta)\cosh(2\rho\sqrt{|\Delta|})+1}\]
hence the desired result holds.

\end{proof}




Thanks to the above estimation, we get a practically useful version
of Theorem \ref{convergence1}.

\begin{theorem}
Let $(a_k)_k$ be a sequence in $(0,1]$ such that
\[\lim_{k\rightarrow\infty}a_k=0\quad\text{and}\quad\sum_{k=0}^{\infty}a_k=+\infty\]
and let $t_k=\beta a_k$ in the algorithm with
\[0<\beta\leq\frac{\rho-\sigma}{C(\rho,\delta)\cosh(2\rho\sqrt{|\Delta|})+1}\]
then we have
\[\lim_{k\rightarrow\infty}d(x_k,\mathfrak{M}_{\mu})=0\quad\text{and}\quad\lim_{k\rightarrow\infty}f(x_k)=f_*\]
Moreover, if $(a_k)_k$ verifies also that
\[\sum_{k=0}^{\infty}a^2_k<+\infty\]
then there existes some $m\in \mathfrak{M}_{\mu}$ such that
$x_k\rightarrow m$, when $k\rightarrow\infty$.
\end{theorem}

\begin{proof}
This is a simple corollary to Theorem \ref{convergence1} and Lemma
\ref{estimation of rx}.
\end{proof}

Now we turn to the problem of error estimates of the subgradient
algorithm under condition $\ast$.





\begin{proposition}\label{error}

Let condition $\ast$ hold and the stepsize $(t_k)_k$ satisfy
 \[\lim_{k\rightarrow\infty}t_k=0\quad \text{and}\quad\sum_{k=0}^{\infty}t_k=+\infty\]
Then there exists $N\in\mathbf{N}$, such that for every $k\geq N$,
\[d^2(x_k,m)\leq b_k\]
where $m$ is the unique median of $\mu$ and the sequence
$(b_k)_{k\geq N}$ is defined by
\[b_N=(\rho+\sigma)^2\quad\text{and}\quad b_{k+1}=(1-2\tau
t_k)b_k+C(\rho,\delta)t^2_k\,,\quad k\geq N\] which converges to $0$
when $k\rightarrow\infty$. More explicitely, for every $k\geq N$,
\[b_{k+1}=(\rho+\sigma)^2\prod_{i=N}^{k}(1-2\tau t_i)
+C(\rho,\delta)\big(\sum_{j=N+1}^{k}t^2_{j-1}\prod_{i=j}^{k}(1-2\tau
t_i)+t^2_k\,\,\big)\]
\end{proposition}

\begin{proof}
Since $t_k\rightarrow0$, there exists $N\in\mathbf{N}$, such that
for every $k\geq N$, we have $2\tau t_k<1$. For every of these $k$,
we choose a geodesic $\gamma$: $[0,1]\rightarrow\bar{B}(a,\rho)$
with $\gamma(0)=m$ and $\gamma(1)=x_k$. Then by Theorem \ref{strong
convex},
\[f(x_k)-f_*\geq
\frac{d}{ds}f(\gamma(s))\big|_{s=0+}+\tau d^2(x_k, m)\geq \tau
d^2(x_k, m)\] combining this and Lemma \ref{fundamental lemma} we
get
\[d^2(x_{k+1},m)\leq(1-2\tau t_k)d^2(x_k,m)+C(\rho,\delta)t^2_k\]
Proposition \ref{location} yields that $d^2(x_N,m)\leq
(\rho+\sigma)^2=b_N$ and then by induction it is easily seen that
$d^2(x_k,m)\leq b_k$ for every $k\geq N$. The same method as in the
proof of Theorem \ref{convergence1} can show that $b_k\rightarrow0$.
In fact, we first show that
\[\liminf_{k\rightarrow\infty} b_k=0\]
If this is not true, then there exist $N_1\geq N$ and $\eta>0$ such
that for every $k\geq N_1$ we have $b_k>\eta$. Since
$\lim_{k\rightarrow\infty}t_k=0$, we can suppose that
$C(\rho,\delta)t_k<\tau\eta$ for every $k\geq N_1$ and then
\[b_{k+1}=b_k+t_k(C(\rho,\delta)t_k-2\tau b_k)\leq b_k+t_k(C(\rho,\delta)t_k-2\tau\eta)\leq b_k-\tau\eta t_k\]
by summing the above inequalities we get
\[\tau\eta\sum_{i=N_1}^{k}t_i\leq b_{N_1}-b_{k+1}\leq b_{N_1}\]
which contradicts with $\sum_{k=0}^{\infty}t_k=+\infty$, this proves the assertion.\\
Then for every $k\geq N$, we consider the following two cases:\\
If $b_k>C(\rho,\delta)t_k/(2\tau)$, then we get
\[b_{k+1}\leq b_k-2\tau
t_k(C(\rho,\delta)t_k/(2\tau))+C(\rho,\delta)t^2_k=b_k\] If $b_k\leq
C(\rho,\delta)t_k/(2\tau)$, then we get
\[b_{k+1}=(1-2\tau
t_k)C(\rho,\delta)t_k/(2\tau)+C(\rho,\delta)t^2_k=C(\rho,\delta)t_k/(2\tau)\]
Thus we have always that
\[b_{k+1}\leq\max\{b_k,C(\rho,\delta)t_k/(2\tau)\}\]
by induction we get for every $n\geq k$,
\[b_{n+1}\leq
\max\{b_k,(C(\rho,\delta)/(2\tau))\max\{t_k,t_{k+1},\dots,t_n\}\}\]
thus we have
\[\limsup_{n\rightarrow\infty}b_{n}\leq \max\{b_k,(C(\rho,\delta)/(2\tau))\sup\{t_i:i\geq k\}\}\]
In order to get $b_k\rightarrow0$, it suffices to take the inferior
limite on the right hand side. Finally, the explicit expression of
$(b_k)_k$ follows from induction.
\end{proof}

We proceed to show that if the stepsize $(t_k)_k$ is chosen to be
the harmonic series, then the rate of convergence of our algorithm
is sublinear. To do this, we use the following Lemma in Nedic and
Bertsekas \cite{Nedic1}.

\begin{lemma}
Let $(u_k)_k$ be a sequence of nonnegative real numbers satisfying
\[u_{k+1}\leq(1-\frac{\alpha}{k+1})u_k+\frac{\zeta}{(k+1)^2},\quad k\geq0\]
where $\alpha$ and $\zeta$ are positive constants. Then
\[
u_{k+1}\leq
\begin{cases}

\frac{1}{(k+2)^\alpha}(u_0+\frac{2^{\alpha}\zeta(2-\alpha)}{1-\alpha})
& \text{if $0<\alpha<1$;}\\\\

\frac{\zeta(1+\ln(k+1))}{k+1}  & \text{if $\alpha=1$;}\\\\

\frac{1}{(\alpha-1)(k+2)}(\zeta+\frac{(\alpha-1)u_0-\zeta}{(k+2)^{\alpha-1}})
& \text{if $\alpha>1$;}

\end{cases}
\]
\end{lemma}

\begin{proposition}
Let condition $\ast$ hold and we choose $t_k=r/(k+1)$ for every
$k\geq0$ with some constant $r>0$, then we have
\[
d^2(x_{k+1},m)\leq
\begin{cases}

\frac{1}{(k+2)^\alpha}((\rho+\sigma)^2+\frac{2^{\alpha}r^2C(\rho,\delta)(2-\alpha)}{1-\alpha})
& \text{if
$0<\alpha<1$;}\\\\

\frac{r^2C(\rho,\delta)}{k+1}(1+\ln(k+1)) & \text{if $\alpha=1$;}\\\\

\frac{1}{(\alpha-1)(k+2)}(r^2C(\rho,\delta)
+\frac{(\alpha-1)(\rho+\sigma)^2-r^2C(\rho,\delta)}{(k+2)^{\alpha-1}})
& \text{if $\alpha>1$;}

\end{cases}
\]
where $m$ is the unique median of $\mu$ and $\alpha=2\tau r$
\end{proposition}

\begin{proof}
As in the proof of Proposition \ref{error}, we get that for every
$k\geq0$,
\[d^2(x_{k+1},m)\leq(1-\frac{2\tau r}{k+1})d^2(x_k,m)+\frac{r^2C(\rho,\delta)}{(k+1)^2}\]
then it suffices to use the preceding lemma with $\alpha=2\tau r$
and $\zeta=r^2C(\rho,\delta)$. Observe also that $d(x_0,m)\leq
\rho+\sigma$ by Proposition \ref{location}.
\end{proof}


\begin{thebibliography}{30}

\bibitem{Arnaudon1}
Arnaudon M. Esp\'{e}rances conditionnelles et C-martingales dans les
vari\'{e}t\'{e}s. \emph{S\'{e}minaire de Probabilit\'{e}s}-XXVIII.
300-311. Lecture Notes in Maths. 1583. Springer. Berlin. 1994


\bibitem{Arnaudon2}
Arnaudon M. Barycentres convexes et approximations des martingales
dans les vari\'{e}t\'{e}s. \emph{S\'{e}minaire de
Probabilit\'{e}s}-XXIX. 70-85. Lecture Notes in Maths. 1613.
Springer. Berlin. 1995


\bibitem{Arnaudon3}
Arnaudon M. and Li X.M. Barycenters of measures transported by
stochastic flows. \emph{Annals of probability}. vol 33. no. 4. 2005.
pp 1509-1543.


\bibitem{Barbaresco1}
Barbaresco F. Innovative Tools for Radar Signal Processing Based on
Cartan's Geometry of SPD Matrices and Information Geometry. IEEE
International Radar Conference (2008).


\bibitem{Barbaresco2}
Barbaresco F. Interactions between Symmetric Cone and Information
Geometries, ETVC'08, Springer Lecture Notes in Computer Science,
vol.5416, pp. 124-163, 2009


\bibitem{Barbaresco3}
Barbaresco F. and Bouyt G. Espace Riemannien sym¨¦trique et
g\'{e}om\'{e}trie des espaces de matrices de covariance :
\'{e}quations de diffusion et calculs de m\'{e}dianes, GRETSI'09
conference, Dijon, September 2009

\bibitem{Barbaresco4}
Barbaresco F. New Foundation of Radar Doppler Signal Processing
based on Advanced Differential Geometry of Symmetric Spaces: Doppler
Matrix CFAR and Radar Application, Radar'09 Conference, Bordeaux,
October 2009

\bibitem{Correa}
Correa R. and Lemar\'{e}chal C. Convergence of some algorithms for
convex minimization. Mathematical Programming 62 (1993) 261-275
North-Holland.


\bibitem{Emery}
Emery M. Mokobodzki G. Sur le barycentre d'une probabilit\'{e} dans
une vari\'{e}t\'{e}. \emph{S\'{e}minaire de Probabilit\'{e}s
(Strasbourg)} tome 25 (1991) p.220-233.


\bibitem{Ferreira}
Ferreira O.P. and Oliveira P.R.  Subgradient Algorithm on Riemannian
Manifolds. Journal of Optimization Theory and Applications: Vol 97,
No. 1, pp. 93-104, April 1998


\bibitem{Fletcher1}
P. T. Fletcher and S. Joshi. Principle geodesic analysis on
symmetric spaces: statistics of diffusion tensors. Proceedings of
ECCV Workshop on Computer Vision Approaches to Medical Image
Analysis.


\bibitem{Fletcher2}
P. T. Fletcher, C. Lu and S. Joshi. Statistics of shape via
principle geodesic analysis on Lie groups. Proceedings of the IEEE
conference on Computer Vision and Pattern Recognition.


\bibitem{Fletcher3}
P. T. Fletcher, S. Venkatasubramanian and S. Joshi. The geometric
median on Riemannian manifolds with application to robust atlas
estimation. NeuroImage 45 (2009) S143-S152.


\bibitem{Karcher}
Karcher H. Riemannian center of mass and mollifier smoothing.
Communications on Pure and Applied Mathematics. vol xxx. 509-541
(1977)


\bibitem{Kendall}
Kendall W. S. Probability, convexity, and harmonic maps with small
image I: uniqueness and fine existence. Proc. London Math. Soc. (3)
61 (1990) 371-406


\bibitem{Khun}
Kuhn. H.W. A note on Fermat's problem. Mathematical Programming 4
(1973) 98-107 North-Holland Publishing Company.


\bibitem{Le}
Le H. Estimation of Riemannian barycenters. LMS J. Comput. Math. 7
(2004) 193-200


\bibitem{Nedic1}
Angelia Nedic and Dimitri P. Bertsekas.  Convergence Rate of
Incremental Subgradient Algorithms. Stochastic Optimization:
Algorithms and Applications(S.Uryasev and P. M. Pardalos, Editors),
pp. 263-304. Kluwer Academic Publishers 2000.


\bibitem{Nedic2} Angelia Nedic and Dimitri P. Bertsekas.
Incremental Subgradient Methods for Nondifferentiable Optimization.
SIAM J. Optim. Volume 12, Issue 1, pp. 109-138(2001)


\bibitem{Ostresh}
Ostresh L.M. JR. On the convergence of a class of iterative methods
for solving Weber location problem. Operation Research. vol 26. No.
4, July-August 1978.


\bibitem{Pennec}
Pennec X. Intrinsic statistics on Riemannian manifolds: Basic tools
for geometric measurements. \emph{Journal of Mathematical Imaging
and Vision} 2006


\bibitem{Sahib}
Sahib A. Esp\'{e}rance d'une variable al\'{e}atoire \`{a} valeur
dans un espace m\'{e}trique. Th\`{e}se de l'Universit\'{e} de Rouen
(1998).



\bibitem{Weiszfeld}
Weiszfeld E. Sur le point pour lequel la somme des distances de n
points donn\'{e}s est minimum. Tohoku Math. J. 43 (1937), 355-386.


\end{thebibliography}
\end{document}